\documentclass[11pt]{amsart}
\usepackage[hmargin=2.5cm,vmargin=2.5cm]{geometry}
\usepackage{amsfonts}
\usepackage{amsthm}
\usepackage[english]{algorithm2e}
\usepackage{amsmath}
\usepackage{amssymb}
\usepackage{hyperref}
\usepackage[T1]{fontenc}
\newtheorem{Theorem}{Theorem}[section]
\newtheorem{Definition}[Theorem]{Definition}
\newtheorem{Remark}[Theorem]{Remark}

\newtheorem{Proposition}[Theorem]{Proposition}
\newtheorem{Corollary}[Theorem]{Corollary}

\title{On Hom-Leibniz and Hom-Lie-Yamaguti Superalgebras}
\author[ Donatien Gaparayi, Sylvain Attan and A. Nourou Issa]
       {Donatien Gaparayi $^1,$  Sylvain Attan $^2$ and A. Nourou Issa $^3$
       \\\\
       $^1$ Ecole Normale Sup\'erieure (E.N.S),\\
        BP 6983 Bujumbura, Burundi\\
          $^{2,3}$ D\'epartement de Math\'ematiques\\
       Universit\'e d'Abomey-calavi\\
       01 BP 4521, Cotonou 01, B\'enin\\
        $^1$ gapadona@yahoo.fr, $^2$ syltane2010@yahoo.fr, $^3$ woraniss@yahoo.fr
       }
\begin{document}
\maketitle

\begin{abstract}
 In this paper some characterizations of Hom-Leibniz superalgebras are given and some of their basic properties are found. These properties can be seen as a generalization of corresponding well-known properties of Hom-Leibniz algebras. Considering the Hom-Akivis superalgebra asso-\\ciated
to a given Hom-Leibniz superalgebra, it is observed that the Hom-super Akivis identity leads to an additional property of Hom-Leibniz superalgebras, which in turn gives a necessary and sufficient condition for Hom-super Lie admissibility of Hom-Leibniz superalgebras. We show 
also that every (left) Hom-Leibniz superalgebra has a natural Hom-Lie-Yamaguti superalgebra structure.
\\

AMS {\it Subject Class. (2010): 17A30, 17A32, 17D99}\vspace{0.30cm}
\end{abstract}
\section{Introduction} 

A (left) Leibniz superalgebra is an superalgebra $(L = L_0 \oplus L_1,\ast)$ satisfying the identity
\begin{equation}
x\ast(y\ast z) = (x\ast y)\ast z +(-1)^{|x||y|} y\ast(x\ast z).
\end{equation}
As Leibniz algebras which are been  introduced by J. -L. Loday \cite{Loday} (and so they are sometimes called Loday algebras) as a noncommutative analogue of Lie algebras,  Leibniz superalgebras are also seen as  supernoncommutative analogues of  Lie 
superalgebras. Indeed, if the operation $"\ast"$ of a given Leibniz superalgebra $(L,\ast)$ is skew-supersymmetric, then $(L,\ast)$ is a Lie  superalgebra.

One of the problems in the general theory of a given class of (binary or binary-ternary) nonassociative algebras is the study of relationships
between that class of algebras and the one of Lie algebras. In the same rule, the search of relationships between a class of nonassociative 
algebras and the one of Leibniz algebras is of interest (at least for constructing concrete examples of the given class of nonassociative 
algebras). In this setting, the existence of a Lie-Yamaguti structure on any (left) Leibniz algebra pointed out in \cite{KIN1} is a good 
illustration. The counterpart of this construction in the Hom-algebra setting has investigated in \cite{GAP2}. Indeed, the authors show that
every multiplicative left Hom-Leibniz algebra has a natural Hom-Lie-Yamaguti structure. Hom-Lie-Yamaguti algebras (Hom-LY 
algebras) is introduced in \cite{GAP}. Recall that a (left) Hom-Leibniz algebra \cite{Makhlouf3} is a Hom-algebra $(L,.,\alpha)$ satisfying
the (left) Hom-Leibniz identity
\begin{equation}
\alpha(x)\cdot(y\cdot z) = (x\cdot y)\cdot\alpha(z) + \alpha(y)\cdot (x\cdot z)  
\end{equation}
for all $x,y,z$ in $L$ and Hom-Lie Yamaguti algebra \cite{GAP} (Hom-LY algebra for short) is a quadruple $(L,\ast,\{,,\},\alpha )$  in which $L$ is
$\Bbb K$-vector space, $"\ast"$ a binary operation and $"\{,,\}"$ a ternary operation on $L$, and $\alpha : L\rightarrow L$ a linear map
such that 
\begin{eqnarray*}
&& (HLY1) \,\,\alpha (x \ast y) = \alpha(x)\ast \alpha(y),\\
&& (HLY2) \,\,\alpha (\{x,y,z\}) = \{\alpha(x),\alpha(y),\alpha(z)\},\\ 
&& (HLY3)\,\, x\ast y = - y\ast x,\\ 
&& (HLY4) \,\,\{x,y,z\} = - \{y,x,z\},\\
&& (HLY5) \,\,\circlearrowleft_{(x,y,z)} [(x\ast y)\ast \alpha(z) + \{x,y,z\}] = 0,\\ 
&& (HLY6)\,\, \circlearrowleft_{(x,y,z)}[\{x\ast y, \alpha(z),\alpha(u)\}] = 0,\\
&& (HLY7) \,\,\{\alpha(x),\alpha(y),u\ast v\} = \{x,y,u\} \ast \alpha^2(v) + \alpha^2(u)\ast \{x,y,v\},\\ 
&& (HLY8)\,\, \{\alpha^2(x),\alpha^2(y),\{u,v,w\}\} =  \{\{x,y,u\},\alpha^2(v),\alpha^2(w)\} + \{\alpha^2(u),\{x,y,v\},\alpha^2(w)\}\\
&& \hspace{6.5cm}+ \,\, \{\alpha^2(u),\alpha^2(v),\{x,y,w\}\}, \,\,\forall \,\,u,v,w,x,y,z \in L.
\end{eqnarray*}
Observe that if $x\ast y = 0$, for all $x,y$ in $L$, then Hom-Lie-Yamaguti algebras $(L,\ast,\{,,\},\alpha )$ reduce to a Hom-Lie triple
system $(L, \{, , \},\alpha)$ as defined in \cite{YAU5} and if $\{x,y,z\} = 0$ for all $x, y, z$ in $L$, then $(L,\ast,\{,,\},\alpha )$
is a Hom-Lie algebra $(L,\ast,\alpha)$ \cite{Hartwig}. 

The Hom-algebra structures arise first in quasi-deformation of Lie algebras of vector fields. Discrete modifications of vector fields via 
twisted derivations lead to Hom-Lie and quasi-Hom-Lie structures in which the Jacobi condition is twisted \cite{Hartwig}. A general study and
construction of Hom-Lie algebras have been considered (see \cite{Jin},\cite{Schestakov},\cite{Yuan}). Many researchers
investigated Hom-structures on different objects. Yau introduced Hom-Malcev algebras and gave the connections between Hom-alternative
algebras and Hom-Malcev algebras \cite{YAU1}. For further more information on other Hom-type algebras, one may refer to, e.g.,
\cite{Ataguema1}, \cite{Gohr2}, \cite{Makhlouf1}, \cite{YAU3}, \cite{YAU4}, \cite{YAU5}.

In this paper, after considering some characterizations and basic properties of Hom-Leibniz superalgebras, we show that every (left) 
Hom-Leibniz superalgebra has a natural Hom- Lie-Yamaguti superalgebra structure. We observe that the even part of this investigation recovers to the 
works done in \cite{GAP2} and in \cite{NOU2}.

The rest of this paper is organized as follows. In Section two, we recall basic definitions in the
Hom-superalgebras theory and useful 
 results  about Hom-associative superalgebras and Hom-Lie superalgebras. In \cite{MAK5}, the authors show that the supercommutator bracket
 defined using the multiplication in a Hom-associative superalgebra leads naturally to Hom-Lie superalgebra. In particular, we recall the
 notion of Hom-Leibniz, Hom-Akivis and Hom-Lie-Yamaguti superalgebras. In fact, in Proposition \ref{As-Ak} we show that The 
 supercommutator-Hom-associator superalgebra of a multiplicative non-Hom-associative superalgebra is a Hom-Akivis superalgebra. In the third Section, super versions of some well-known properties of (left) Hom-Leibniz algebras \cite{NOU2} are
 displayed. Consi-\\dering the  specific properties  of the binary and ternary operations of the Hom-Akivis superalgebra associated to a
 given  Hom-Leibniz superalgebra, we infer a characteristic  property of Hom-Leibniz superalgebras. This property in turn allows to give a 
 necessary and sufficient condition  for the super Hom-Lie  admissibility  of these Hom-superalgebras. In the last  Section,  we prove the
 existence of a Hom-Lie-Yamaguti superalgebra structure on any (multiplicative) left Hom-Leibniz superalgebra. Indeed, every multiplicative (left) 
 Hom-Leibniz superalgebra has a natural Hom-Lie-Yamaguti superalgebra structure.
\section{Preliminaries} 
In this section, we recall useful definitions and results that are for further use. In particular, we recall the notion of Hom-Leibniz, 
Hom-Akivis and Hom-Lie-Yamaguti superalgebras.
\begin{Definition}
\begin{itemize}
\item  [$(i)$] Let $f:(A, \ast, \alpha) \rightarrow (A' , \ast',\alpha')$ be a linear map, where $A = A_0\oplus A_1$ and $A' = A'_0\oplus A'_1$ are 
$\mathbf{Z}_2-$graded vector spaces ( $\mathbb{K}$-vectors superspaces). The map $f$ is called an even (resp. odd ) map if  $f(A_i) \subset A'_i$ (resp. $f (A_i ) \subset A'_{i+1})$,
for $i = 0,1.$
\item [$(ii)$] A (multiplicative) \textbf{ $n$-ary Hom-superalgebra} is a triple $(A,\{,\cdots,\}, \alpha)$ in which $A = A_0\oplus A_1$ is a $\Bbb K$-vectors superspace, 
$\{,\cdots,\}: A^n \rightarrow A$ is an $n$-linear map such that $\{A_{i_1}, A_{i_2},\cdots,A_{i_n}\}\subseteq A_{i_1+i_2+\cdots+i_n}$ and $\alpha: A \rightarrow A$ is an even linear map such that \\
$\alpha(\{x_1,x_2,\cdots, x_n \}) = \{\alpha(x_1),\alpha(x_2),\cdots,\alpha(x_n)\}$ (multiplicativity).
\end{itemize}
\end{Definition}
We will be  interested in binary ($n=2$), ternary ($n=3$) and binary-ternary  Hom-algebras (i.e. Hom-algebras with binary and ternary operations). For convenience, throughout this paper  we assume that all (binary, ternary and binary-ternary) Hom-superalgebras are multiplicative and for a given vector superspace $A$, the set of all homogeneous elements will be denoted by $\mathcal{H}(A).$ 
\begin{Definition} \label{H-SAs}
  Let $(A,\ast, \alpha)$ be a Hom-superalgebra.
  \begin{itemize}
 \item [$(i)$] The Hom-associator of $A$ \cite{MAK4} is the trilinear map $as_{\alpha}: A \times A \times A \rightarrow A$ defined as
\begin{equation}\label{asA}
as_{\alpha} = \ast \circ (\ast \otimes \alpha - \alpha \otimes \ast).
\end{equation}
In terms of homogeneous elements, the map $as_{\alpha}$ is given by
\begin{equation*}
as_{\alpha}(x,y,z) = (x\ast y)\ast \alpha(z)- \alpha(x)\ast(y\ast z).
\end{equation*}
\item [$(ii)$]  A Hom-associative superalgebra \cite{MAK5} is a multiplicative Hom-superalgebra $(A,\ast,\alpha)$ such that $$as_{\alpha}(x,y,z) = 0,\,\,\, 
\forall\,\,x,y,z\in \mathcal{H}(A).$$
However, the Hom-associativity is not always assumed, i.e. $as(x,y,z) \neq 0$ in general. In this case $(A,\cdot,\alpha)$ is said 
non-Hom-associative \cite{Issa2} (or Hom-nonassociative \cite{YAU3}; in \cite{MAK2}, $(A,\cdot, \alpha)$ is also called a nonassociative 
Hom-superalgebra). This matches the generalization  of associative superalgebras by the nonassociative ones.
\item [$(iii)$]  The Hom-super-Jacobian of $A$ \cite{MAK4} is the trilinear map $J_{\alpha} : A\times A\times A \rightarrow A$ defined
as $J_\alpha(x,y,z):= (x\ast y)\ast \alpha(z) + (-1)^{|x|(|y|+|z|)}(y\ast z)\ast \alpha(x) + (-1)^{|z|(|x|+|y|)}(z\ast x)\ast \alpha(y)$
for all homogeneous elements $ x,y,z \in A.$ If $\alpha = Id$, then the Hom-super-Jacobian reduces to the usual super-Jacobian.
\item [$(iv)$]  A Hom-Lie superalgebra \cite{MAK5} is a multiplicative Hom-superalgebra $(A,\ast,\alpha)$ such that
\begin{eqnarray}\label{HSJ}\nonumber
&& x\ast y = - (-1)^{|x||y|}y\ast x\\
&& J_{\alpha} (x,y,z) = 0 \,\,\,\mbox{(the Hom-super-Jacobi identity)}
\end{eqnarray}
for all homogeneous elements $x,y,z\in A.$ 
If $\alpha = Id,$ a Hom-Lie superalgebra reduces to an usual Lie superalgebra.
\end{itemize}
\end{Definition}
\begin{Definition} \cite{w1}
 A (left) Hom-Leibniz superalgebra is a multiplicative Hom-superalgebra $(L,\ast,\alpha)$ satisfying the (left) Hom-Leibniz superalgebra identity
 \begin{equation}\label{LLSI}
  \alpha(x)\ast(y\ast z) = (x\ast y)\alpha(z) + (-1)^{|x||y|}\alpha(y)\ast(x\ast z)
 \end{equation}
for all homogeneous elements $x,y,z \in L.$
\end{Definition}
\begin{Remark}
 The dual superidentity of (\ref{LLSI}) is given by
 \begin{equation}\label{RLSI}
  (x\cdot y)\cdot\alpha(z) = \alpha(x)\cdot(y\cdot z) + (-1)^{|y||z|}(x\cdot z)\cdot\alpha(y)
 \end{equation}
which defines a (right) Hom-Leibniz superalgebra. Indeed, if the operation $"\ast"$ verifies (\ref{LLSI}), then the operation $"\cdot"$ 
defined by $x\cdot y = y\ast x$ verifies also (\ref{RLSI}) and then from a (left) Hom-Leibniz superalgebra, one can get the (right) Hom-Leibniz superalgebra.
\end{Remark}
In terms of the Hom-associator, the identity (\ref{LLSI}) is written as
\begin{equation}\label{Ass-HSLei}
as(x, y, z) = - (-1)^{|x||y|}\alpha(y)\ast (x\ast z)
\end{equation}
Therefore, from Definition \ref{H-SAs} \ in assertion $(ii)$, we see that Hom-Leibniz superalgebras are examples of non-Hom-associative superalgebras.
\begin{Definition}\cite{gi}\label{AKI-Super}
 A Hom-Akivis superalgebra is a multiplicative binary-ternary Hom-superalgebra $(A,[,],\{,,\},\alpha)$ such that 
\begin{eqnarray}\label{Sakivs}
 [x,y] &=& (-1)^{|x||y|}[y,x]\nonumber\\
J_{\alpha} (x,y,z)  &=& \circlearrowleft_{(x,y,z)}(-1)^{|x||z|}[x,y,z] 
- \circlearrowleft_{(x,y,z)}(-1)^{(|y|+|z|)|x|}\{y,x,z\} ,\label{SAkivs}
\end{eqnarray}
for all homogeneous elements $x,y,z$ in $A$. 
The identity (\ref{Sakivs}) is called the Hom-super Akivis identity. 
\end{Definition}
Observe that when $\alpha = Id,$ the Hom-super Akivis identity (\ref{Sakivs}) reduces to the usual super Akivis identity (see in \cite{Santana}).
\newline
\begin{Proposition}\label{As-Ak}
 The supercommutator-Hom-associator superalgebra of a multiplicative non-Hom-associative superalgebra is a  Hom-Akivis superalgebra.
\end{Proposition}
 \textbf{Proof.} Let $(A,\ast,\alpha)$ be a multiplicative non-Hom-associative superalgebra. For any homogeneous elements $x,y,z \in A$, define the operations
 
$[x,y] := x\ast y - (-1)^{|x||y|} y\ast x$ (supercommutator)

$\{x,y,z\}_{\alpha} := as_{\alpha}(x,y,z)$ (Hom-associator; see (\ref{asA})).\\
Then

$[[x, y], \alpha(z)] \\= (x\ast y)\ast \alpha(z)- (-1)^{|x||y|}(y\ast x)\ast \alpha(z)- (-1)^{|z|(|x|+|y|)}\alpha(z)\ast(x\ast y)
+ (-1)^{|x|(|z|+|y|)+|y||z|}\alpha(z)\ast(y\ast x)$
 and
 
$\{x,y,z\}_{\alpha} - (-1)^{|x||y|}\{y,x,z\}_{\alpha}\\ = (x\ast y)\ast \alpha(z) - \alpha(x)\ast(y\ast z) 
- (-1)^{|x||y|}(y\ast x)\ast \alpha(z) + (-1)^{|x||y|}\alpha(y)\ast(x\ast z).$
 Expanding $\circlearrowleft_{(x,y,z)}[[x, y],\alpha(z)]$ and $\circlearrowleft_{(x,y,z)}(\{x, y, z\}_{\alpha} 
 - (-1)^{|x||y|}\{y, x, z\}_{\alpha})$ respectively, where $\circlearrowleft_{(x,y,z)}$ denotes the sum over cyclic permutation of $x,y,z,$
 one gets  (\ref{Sakivs}) and so $(A,[,],\{,,\}_{\alpha},\alpha)$ is a Hom-Akivis superalgebra. The multiplicativity of $(A,[,], \{,,\}_{\alpha},\alpha)$ follows from the one of $(A,\ast,\alpha).$ \hfill$\square$ 
 
 Observed that to each non-Hom-associative superalgebra is associated a Hom-Akivis superalgebra (this is the Hom-analogue of a similar 
 relationship between nonassociative superalgebras and Akivis superalgebras \cite{Santana}). In this note, in the third Section
 we use the specific properties of the Hom-Akivis superalgebra associated to a given Hom-Leibniz superalgebra to derive a property
 characterizing Hom-Leibniz superalgebras. Before this, we recall the following.
 
 \begin{Definition}\cite{ga}\label{SHLY-Def}
 {\it A \textbf{Hom-Lie Yamaguti superalgebra} (Hom-LY superalgebra for short) is a binary-ternary Hom-superalgebra $(L,\ast,\{,,\},\alpha )$  such that
\begin{eqnarray*}
&& (SHLY1) \,\,\alpha (x \ast y) = \alpha(x)\ast \alpha(y),\\
&& (SHLY2) \,\,\alpha (\{x,y,z\}) = \{\alpha(x,)\alpha(y),\alpha(z)\},\\ 
&& (SHLY3)\,\, x\ast y = - (-1)^{|x||y|}y\ast x,\\ 
&& (SHLY4) \,\,\{x,y,z\} = - (-1)^{|x||y|}\{y,x,z\},\\
&& (SHLY5) \,\,\circlearrowleft_{(x,y,z)} (-1)^{|x||z|}[(x\ast y)\ast \alpha(z) + \{x,y,z\}] = 0,\\ 
&& (SHLY6)\,\, \circlearrowleft_{(x,y,z)}(-1)^{|x||z|}[\{x\ast y, \alpha(z),\alpha(u)\}] = 0,\\
&& (SHLY7) \,\,\{\alpha(x),\alpha(y),u\ast v\} = \{x,y,u\} \ast \alpha^2(v) + (-1)^{|u|(|x|+|y|)}\alpha^2(u)\ast \{x,y,v\},\\ 
&& (SHLY8)\,\, \{\alpha^2(x),\alpha^2(y),\{u,v,w\}\} =  \{\{x,y,u\},\alpha^2(v),\alpha^2(w)\}\\
 && \hspace{6.5cm}+\,\,  (-1)^{|u|(|x|+|y|)}\{\alpha^2(u),\{x,y,v\},\alpha^2(w)\}  \\
 && \hspace{6.5cm}+ \,\, (-1)^{(|x|+|y|)(|u|+|v|)}\{\alpha^2(u),\alpha^2(v),\{x,y,w\}\},
\end{eqnarray*}
for all homogeneous elements $u,v,w,x,y,z \in L$ and where $\circlearrowleft_{(x,y,z)}$ denotes the sum over cyclic permutation of $x,y,z.$ } 
\end{Definition}
Note that the conditions {$(SHLY1)$} and {$(SHLY2)$} mean the multiplicativity of $(L,\ast,\{,,\},\alpha ).$
 
 \section{Some characterizations of Hom-Leibniz superalgebras}
 In this section, super-versions of some well-known properties of (left) Hom-Leibniz algebras \cite{NOU2} are displayed. Considering the 
 specific properties  of the binary and ternary operations of the Hom-Akivis superalgebra associated to a given Hom-Leibniz superalgebra,
 we infer a characteristic  property of Hom-Leibniz superalgebras. This property in turn allows to give a necessary and sufficient condition
 for the Hom-super Lie  admissibility  of these Hom-superalgebras.
 
 \begin{Proposition}\label{H-SLei-SAkiv}
 Let $(L, \ast, \alpha)$ be a  Hom-Leibniz superalgebra. Consider on $(L, \ast, \alpha)$ the operations
 \begin{eqnarray}
 && [x,y] := x \ast y - (-1)^{|x||y|}y \ast x\\
&& \{x, y, z\} := as(x, y, z). 
 \end{eqnarray}
Then $(L,[,],\{,,\})$ is a Hom-Akivis superalgebra.
\end{Proposition}

\textbf{Proof.} The proof follows if observe that Proposition \ref{H-SLei-SAkiv} is Proposition \ref{As-Ak} using (\ref{Ass-HSLei}).
Then, we have 
 \begin{equation}\label{SAkid-SLei}
 J_{\alpha}(x, y, z) = \,\,\circlearrowleft_{(x,y,z)} (-1)^{|x||z|}(x \ast y) \ast \alpha(z),
 \end{equation}
for all homogeneous elements $x, y, z \in A.$ \hfill$\square$

One observes that (\ref{SAkid-SLei}) is the specific form of the Hom-super Akivis identity (\ref{SAkivs}) in case of Hom-Leibniz superalgebras.
The relation above will be used in the Section 4. The  skew-supersymmetry of the operation $"\ast"$ of a Hom-Leibniz superalgebra $(L,\ast, \alpha)$
is a condition for $(L,\ast, \alpha)$ to be a Hom-Lie superalgebra \cite{MAK5}. From Proposition \ref{H-SLei-SAkiv} one gets the following 
necessary and sufficient condition for the Hom-super Lie admissibility \cite{MAK5} of a given Hom-Lie superalgebra. 

\begin{Corollary}
 A Hom-Leibniz superalgebra $(L,\ast, \alpha)$ is Hom-super Lie admissible if and only if 
 $\circlearrowleft_{(x,y,z)} (-1)^{|x||z|}(x \ast y) \ast \alpha(z) = 0$, for all homogeneous elements $x, y, z$ in $L$. 
\end{Corollary}

In the following, some  properties of Hom-Leibniz superalgebras are found.
\begin{Proposition}\label{PropLLS}
Let $(L,\ast, \alpha)$ be a Hom-Leibniz superalgebra. Then 
\begin{eqnarray}\label{Prop1}
&& (i)\,\,\,(x \ast y + (-1)^{|x||y|}y \ast x) \ast \alpha(z) = 0,\\\label{Prop2}
&& (ii)\,\,\,\alpha(x) \ast [y, z] = [x \ast y, \alpha(z)] + (-1)^{|x||y|}[\alpha(y), x \ast z],
\end{eqnarray}
for all homogeneous elements $x, y, z \in L.$
\end{Proposition}
\textbf{Proof.} The identity (\ref{LLSI}) implies that $(x\ast y) \ast \alpha(z) = \alpha(x) \ast (y\ast z) - (-1)^{|x||y|} \alpha(y) 
\ast (x\ast z).$ Likewise, interchanging $x$ and $y$, we have $(y\ast x) \ast \alpha(z) = \alpha(y) \ast (x\ast z) - \alpha(x) \ast (y\ast z).$
Then, adding memberwise the equalities above, we come to the property $(i)$. Next we have
\begin{eqnarray*}
[x \ast y, \alpha(z)] + (-1)^{|x||y|}[\alpha(y), x \ast z] &=& (x\ast y)\ast \alpha(z) - (-1)^{|z|(|x|+|y|)}\alpha(z) \ast(x\ast y) \\
&+& (-1)^{|x||y|}\alpha(y) \ast (x\ast z)-(-1)^{|y||z|} (x\ast z)\ast \alpha(y)\\
&=& \alpha(x) \ast (y\ast z) - (-1)^{|z|(|x|+|y|)}\alpha(z) \ast (x\ast y)\\
&-& (-1)^{|y||z|}(x\ast z\ast \alpha(y)\,\,\, \mbox{(by (\ref{LLSI}))}\\
&=& \alpha(x) \ast(y\ast z) - (-1)^{|z|(|x|+|y|)}(z\ast x) \ast \alpha(y)\\
&+& (-1)^{|x||z|}\alpha(x) \ast (z\ast y) - (-1)^{|y||z|}(x\ast z) \ast \alpha(y) \,\,\,\mbox{(by (\ref{LLSI}))}\\
&=& \alpha(x) \ast(y\ast z) - (-1)^{|y||z|}\alpha(x) \ast (z\ast y)\,\,\, \mbox{(by (i))}\\
&=& \alpha(x) \ast [y, z].
\end{eqnarray*}
and so we get $(ii)$. $ \hfill \square$ 

\begin{Remark}
 In Proposition \ref{PropLLS}, if  $x,y,z$  are in $L_0$, then one recovers the well-known properties of Hom-Leibniz algebras 
 (see the Proposition 3.1. in \cite{NOU2}).
 \end{Remark}

\section{Hom-Lie-Yamaguti superalgebra on Hom-Leibniz superalgebra}
In this section, we prove the existence of a Hom-Lie-Yamaguti superalgebra structure on any (multiplicative) left 
Hom-Leibniz superalgebra. The proof is based on a specific ternary operation that can be considered on a given Hom-Leibniz superalgebra (this product generalizes the ternary operation considered in \cite{GAP2} on a left Hom-Leibniz algebra which generalizes the one in \cite{KIN1} on a left Leibniz algebra $L$ that produces, along with the skew-symmetrization, a Lie-Yamaguti algebra structure on $L$). Also note that our proof below essentially relies on some properties characterizing Hom-Leibniz superalgebras,
obtained in Proposition \ref{PropLLS}. In the following we have
\begin{Theorem}
 Every  left Hom-Leibniz superalgebra has a natural Hom-Lie-Yamaguti superalgebra structure.
\end{Theorem}

\textbf{Proof.} In a left Hom-Leibniz algebra $(L,\ast,\alpha)$ consider the
 skew-supersymmetrization $$[x, y] := x \ast y - (-1)^{|x||y|} y\ast x$$
for all homogeneous elements $x, y$ in L. In the following, consider the left translations $\Lambda_a b := a \ast b$ in $(L,\ast,\alpha)$.
Then the identities (\ref{LLSI}) and (\ref{Prop2}) can be written respectively as
\begin{eqnarray}\label{TLSR0}
&&\Lambda_{\alpha(x)} (y \ast z) = (\Lambda_x y) \ast \alpha(z) + (-1)^{|x||y|}\alpha(y) \ast (\Lambda_x z), \\\label{TLSR1}
&&\Lambda_{\alpha(x)}[y,z] = [\Lambda_x y,\alpha(z)] + (-1)^{|x||y|}[\alpha(y), \Lambda_x z]. 
\end{eqnarray}
In $(L,\ast,\alpha)$ consider the following ternary operation:
\begin{equation}\label{TLSR2}
\{x, y, z\} := as_{\alpha} (y, x, z) - as_{\alpha} (x, y, z)
\end{equation}
for all homogeneous elements $x,y,z$ in $L.$ Then (\ref{TLSR2}), (\ref{LLSI}) and (\ref{Ass-HSLei}) imply
\begin{equation}\label{TLSR3}
\{x, y, z\} = -(x \ast y)\ast \alpha(z).
\end{equation}
Moreover, we have
\begin{eqnarray}\nonumber
[x, y]\ast \alpha(z) & = & (x \ast y - (-1)^{|x||y|}y \ast x)\ast \alpha(z)\\\nonumber \label{TLSR4}
& = & 2(x \ast y)\ast \alpha(z) \,\,\mbox{(by (\ref{Prop2}))}\\
& = &-2\{x, y, z\} \,\,\mbox{(see (\ref{TLSR3}))}
\end{eqnarray}
so that
\begin{equation}\label{TLSR5}
\{x, y, z\} = -\frac{1}{2}[x, y]\ast \alpha(z)
\end{equation}
Thus (\ref{TLSR2}), (\ref{TLSR3}) and (\ref{TLSR5}) are different expressions of the operation ''$\{,,\}$'' that are for use in what follows.
Now we proceed to verify the validity on $(L,\ast,\alpha)$ of the set of identities $(SHLY1)-(SHLY8).$
The multiplicativity of $(L,\ast,\alpha)$ implies $(SHLY1)$ and $(SHLY2)$ while $(SHLY3)$ is the skew-supersymmetrization and $(SHLY4)$ clearly 
follows from (\ref{TLSR2}) (or (\ref{TLSR3})). Next, observe that $(HLY5)$ is just the Hom-super Akivis identity (\ref{Sakivs}) for $(L,\ast,\alpha).$

 Consider now $\circlearrowleft_{(x,y,z)}\{[x,y],\alpha(z),\alpha(u)\}$. 
Then 
\begin{eqnarray*}
&&\circlearrowleft_{(x,y,z)}(-1)^{|x||z|}\{[x,y], \alpha(z),\alpha(u)\}\\
& = & \circlearrowleft_{(x,y,z)}(-1)^{|x||z|}(-[x,y]\ast \alpha(z))\ast \alpha^2(u) \quad(\mbox{by (\ref{TLSR3})})\\
& = & 2 (\circlearrowleft_{(x,y,z)}(-1)^{|x||z|}\{x,y,z\})\ast \alpha^2(u) \quad(\mbox{by (\ref{TLSR5})})\\
& = & -2 ((-1)^{|x||z|}(x\ast y)\ast \alpha(z) + (-1)^{|x||y|}(y\ast z)\ast \alpha(x) + (-1)^{|y||z|}(z\ast x)\ast \alpha(y))\ast \alpha^2(u)\\
& = & -2 ((-1)^{|x||z|}\alpha(x)\ast (y\ast z) - (-1)^{|x|(|y|+|z|)}\alpha(y)\ast (x\ast z) + (-1)^{|x||y|}(y\ast z)\ast \alpha(x)\\
& & +(-1)^{|y||z|}(z\ast x)\ast \alpha(y))\ast \alpha^2(u)\quad(\mbox{by (\ref{LLSI})})\\
& = & -2 (-1)^{|x||z|}(\alpha(x)\ast (y\ast z) + (-1)^{|x|(|y|+|z|)}(y\ast z)\ast \alpha(x))\ast \alpha^2(u)\\
&& -2 (- (-1)^{|x|(|y|+|z|)}\alpha(y)\ast (x\ast z) + (-1)^{|y||z|}(z\ast x)\ast \alpha(y))\ast \alpha^2(u)\\
& = & -2 (- (-1)^{|x|(|y|+|z|)}\alpha(y)\ast (x\ast z) + (-1)^{|y||z|}(z\ast x)\ast \alpha(y))\ast \alpha^2(u) \quad(\mbox{by (\ref{Prop2})})\\
& = & -2 (- (-1)^{|x|(|y|+|z|)}\alpha(y)\ast(x\ast z) -(-1)^{|z|(|x|+|y|)} (x\ast z)\ast \alpha(y))\ast \alpha^2(u) \quad(\mbox{by (\ref{Prop2})})\\
& = & 2 (-1)^{|x|(|y|+|z|)}(\alpha(y)\ast (x\ast z) + (-1)^{|y|(|x|+|z|)}(x\ast z)\ast \alpha(y))\ast \alpha^2(u)\\ 
& = & 0 \quad(\mbox{by (\ref{Prop2})})
\end{eqnarray*}
so that we get $(SHLY6).$ Next 
\begin{eqnarray*}
 \{\alpha(x),\alpha(y),[u,v]\} & = & - \alpha(x\ast y)\ast \alpha([u,v])\quad(\mbox{by (\ref{TLSR3}) and multiplicativity})\\ 
& = & \Lambda_{(-\alpha(x\ast y))}[\alpha(u),\alpha(v)] \\
& = &[\Lambda_{(-x\ast y)}\alpha(u),\alpha^2(v)] + (-1)^{|u|(|x|+|y|)}[\alpha^2(u),\Lambda_{(-x\ast y)}\alpha(v)] \quad(\mbox{by (\ref{TLSR0})})\\
& = & [\{x,y,u\},\alpha^2(v)]  + (-1)^{|u|(|x|+|y|)}[\alpha^2(u),\{x,y,v\}]\quad(\mbox{by (\ref{TLSR3})})
\end{eqnarray*}
which is $(SHLY7).$ Finally, we compute
\begin{eqnarray*}
&&\{\{x,y,u\},\alpha^2(v),\alpha^2(w)\} + (-1)^{|u|(|x|+|y|)}\{\alpha^2(u),\{x,y,v\},\alpha^2(w)\}\\&& 
+(-1)^{(|u|+|v|)(|x|+|y|)}\{\alpha^2(u),\alpha^2(v),\{x,y,w\}\}\\
& = & \{-\Lambda_{x\ast y}\alpha(u),\alpha^2(v),\alpha^2(w)\} + (-1)^{|u|(|x|+|y|)}\{\alpha^2(u),-\Lambda_{x\ast y}\alpha(v),\alpha^2(w)\} \\
&& + (-1)^{(|u|+|v|)(|x|+|y|)}\{\alpha^2(u),\alpha^2(v),-\Lambda_{x\ast y}\alpha(w)\}\\
& = & ((-\Lambda_{x\ast y}\alpha(u))\ast \alpha^2(v)).\alpha^3(w) - (-1)^{|u|(|x|+|y|)}(\alpha^2(u)\ast (-\Lambda_{x\ast y}\alpha(v)))\ast \alpha^3(w) \\
&&- (-1)^{(|u|+|v|)(|x|+|y|)}(\alpha^2(u)\ast \alpha^2(v))\alpha(-\Lambda_{x\ast y}\alpha(w))\\
& = & (\Lambda_{\alpha(x\ast y)}\alpha(u\ast v))\ast \alpha^3(w) 
+ (-1)^{(|u|+|v|)(|x|+|y|)}\alpha^2(u\ast v)\ast \Lambda_{\alpha(x\ast y)}\alpha^2(w)\quad(\mbox{by (\ref{TLSR0}) and multiplicativity})\\
& = & \Lambda_{\alpha^2(x\ast y)}(\alpha(u\ast v)\ast \alpha^2(w)) \quad(\mbox{by (\ref{TLSR1})})\\
& = & -\alpha^2(x\ast y)\ast (-\alpha(u\ast v)\ast \alpha^2(w))\\
& = & -(\alpha^2(x)\ast \alpha^2(y))\ast \alpha(-(u\ast v)\ast \alpha(w))\quad(\mbox{by multiplicativity})\\
& = & \{\alpha^2(x),\alpha^2(y),\{u,v,w\}\} \quad(\mbox{by (\ref{TLSR3})})
\end{eqnarray*}
Therefore $(L,[,],\{,,\},\alpha)$ is a Hom-LY superalgebra.
This completes the proof. \hfill$\square$ \newline


\end{document}